\documentclass[12pt,leqno]{article}
\usepackage{amsmath,amssymb,latexsym}
\usepackage{amscd}
\usepackage{mathrsfs}
\usepackage{srcltx}
\usepackage{theorem}
\usepackage{amsfonts}

\newcommand{\R}{\mathbb{R}}        
\newcommand{\N}{\mathbb{N}}        
\newcommand{\Rn}{\mathbb{R}^d}     

\renewcommand{\DH}{\cD\,'_H(\Rn)}
\newcommand{\DHO}{\cD\,'_H(\Om)}
\newcommand{\Dh}{\cD\,'_H(\R)}

\newcommand{\OC}{\cO_C'(\Rn)}
\newcommand{\Oc}{\cO_C'(\R)}
\newcommand{\OH}{\cO_H'(\Rn)}
\newcommand{\Oh}{\cO_H'(\R)}

\newcommand{\Di}{\cD(\Rn)}
\newcommand{\Dip}{\cD\,'(\Rn)}
\newcommand{\NZ}{\Rn_*}
\newcommand{\Dop}{\cD\,'(\Om)}
\newcommand{\Do}{\cD(\Om)}

\newcommand{\Dto}{\cD(\tOm)}

\newcommand{\Proof}{\textbf{Proof:} \ }
\newcommand{\qed}{\hspace*{\fill} $\Box $}

\newcommand{\cD}{\mathscr{D}}

\newcommand{\cE}{\mathscr{E}}

\newcommand{\cO}{\mathscr{O}}  

\newcommand{\be}{{\bf 1}}

\newcommand{\om}{\omega}

\newcommand{\ssubset}{\subset\subset}

\newcommand{\cM}{\mathcal{M}}

\newtheorem{proposition}{Proposition}[section]
\newtheorem{example}[proposition]{Example}
\newtheorem{lemma}[proposition]{Lemma}
\newtheorem{corollary}[proposition]{Corollary}
\newtheorem{theorem}[proposition]{Theorem}
\newtheorem{remark}[proposition]{Remark}

\newtheorem{definition}{Definition}

\newcommand{\Om}{{\Omega}}

\newcommand{\eps}{\varepsilon}
\newcommand{\fa}{{\,\,\,\forall\,\,}}
\newcommand{\ex}{{\,\,\,\exists\,\,}}

\newcommand{\supp}{\mathrm{supp}\, }

\newcommand{\vp}{\varphi}

\newcommand{\tOm}{\widetilde{\Om}}
\newcommand{\tM}{\widetilde{M}}

\title{{\sc Hadamard operators on $\Dop$
}}
\author{Dietmar Vogt}
\date{}

\begin{document}

\maketitle

\footnotetext{
    \hfil\break\indent \begin{minipage}[t]{14cm}{\em Key words and phrases:} Operators on distributions, operators on test functions, monomials as eigenvectors, spaces of distributions. \end{minipage}.
     \hfil\break\indent
{}}

\begin{abstract} For open sets $\Om\subset\Rn$ we study Hadamard operators on $\Dop$, that is, continuous linear  operators which admit all monomials as eigenvectors. We characterize them as operators of the form $L(S)=S\star T$ where $T$ is a distribution and $\star$ the multiplicative convolution. This extends previous results for the case of $\Om=\Rn$ but requires essentially different methods.

\end{abstract}


In the present note we study Hadamard operators on $\Dop$, $\Om$ open in $\Rn$, that is, continuous linear operators on $\Dop$ which admit all monomials as eigenvectors. We continue research begun in \cite{Vdprime} where we gave a complete characterization of the Hadamard operators in $\Dip$ as follows: an operator $L\in L(\Dip)$ is of Hadamard type if and only if $L(S)=S\star T$ where $S\star T$ is the multiplicative convolution (see below) and $T$ is a distribution the support of which has positive distance to all coordinate hyperplanes and which has a certain behaviour in infinity, more precisely, which is in the class $\OH$ introduced in \cite{Vdprime}. The case of general open $\Om$ is not a mere generalization, but requires essentially different methods.

We study Hadamard operators by means of its transpose $M:=L^*$ which is an operator in $\Do$ with certain properties.
We show first that $M$ can be extended to an operator $\tM\in L(\Dto)$, where $\tOm=\bigcup_{a\in\NZ}a\Om$ is the closure of $\tOm$ under invertible dilations. For $\tM$ we show the existence of a distribution $T\in\Dto$, which has positive distance to all coordinate hyperplanes, such that $\tM\vp = T_x\vp(x\,\cdot)$ for all $\vp\in\Dto$. $T$ then has to fulfill further conditions on its support which describe that $\tM$ maps $\Do$ into $\Do$. We obtain that $L(S)=S\star T$ and so a characterization of the Hadamard operators on $\Do$. This characterization we further evaluate in the case of $0\in\Om$, here $\tOm=\R^d$ hence $T\in\OH$, and in the case of $\Om\subset\NZ$, here $T$ must be compact. Note that for $d=1$ those are the only cases.

The study of Hadamard operators, or Hadamard multipliers, has a long history. In recent times they were studied in Doma\'nski-Langenbruch \cite{DLre,DLalg,DLhad}, Doma\'nski-Langenbruch-Vogt \cite{DLV} for operators on spaces of real analytic functions and in Vogt \cite{Vhad,Vconv} for spaces of $C^\infty$-functions. The arguments in both cases are quite similar. The study of Hadamard operators on spaces of distributions requires quite different methods and the results are quite different.

Throughout the paper we will set $xy=(x_1y_1,\dots,x_dy_d)$, that is,  the coordinatewise multiplication. Its unit is the vector $\be=(1,\dots,1)$. A map $x\mapsto xy$ with fixed $y\in\NZ$ is called a (invertible) dilation. We set $\R_+=]0,+\infty)$ and $Q_+=\R_+^d$.
For two distributions $S$ and $T$ we define $S\star T$ by $(S\star T)\vp = S_y(T_x\vp(xy))$ for all $\vp\in\Di$ for which this formula makes sense.

We use standard notation of Functional Analysis, in  particular, of distribution theory. For unexplained notation we refer to \cite{DK}, \cite{MV}, \cite{LSI}, \cite{LS}.

\section{Basics}

Let $\Om\subset\Rn$ be open and non-empty.

\begin{definition}\label{d1} A map $L\in L(\Dop)$ is called a Hadamard operator if it admits all monomials as eigenvectors. The set of Hadamard operators we denote by $\cM(\Om)$.
\end{definition}

Examples of Hadamard operators are the \em Euler operators, \rm these are operators of the form $P(\theta)$ where $P$ is a polynomial and $\theta_j=x_j\frac{\partial}{\partial x_j}$. They are in $\cM(\Om)$ for every open $\Om$ and they are the only ones with this property (see Theorem \ref{t9}).

From \cite{Vdprime} we take the following definitions.

\begin{definition}\label{d3} $T\in \OH$ if for any $k$ there are finitely many functions $t_\beta$ such that $(1+|x|^2)^{k/2}t_\beta\in L_\infty(\Rn)$ and such that $T=\sum_\beta \theta^\beta t_\beta$.
\end{definition}

We compare $\OH$ with the space $\OC$ of L. Schwartz, which may be defined by any of the following equivalent properties (see \cite[\S 5, Th\'eor\`eme IX]{LS}):
\begin{enumerate}
\item For any $k$ there are finitely many functions $t_\beta$ such that $(1+|x|^2)^{k/2}t_\beta\in L_\infty(\Rn)$ and such that $T=\sum_\beta \partial^\beta t_\beta$.
\item For any $\chi\in\Di$, $T*\chi$ is a rapidly decreasing continuous function.
\end{enumerate}

We have $\Oc\subset \Oh$ and $\Oc\neq\Oh$. An example is $T=e^{-ix}$ which is in $\Oh$ but not in $\Oc$ (see \cite[Section 3]{Vdprime}).

\begin{definition}\label{d2} By $\DH$ we denote the set of all $T\in\OH$ the support of which has positive distance to all coordinate hyperplanes.
\end{definition}

From \cite[Section 3]{Vdprime} we cite: If  $T\in\OC$ and the support of $T$ has positive distance to all coordinate hyperplanes, then $T\in\DH$.

The main result of \cite{Vdprime} is the following theorem. Here $\sigma(x)=\prod_j \frac{x_j}{|x_j|}$ for $x\in\NZ$.

\begin{theorem} $L\in\cM(\Rn)$ if and only if there is a distribution $T\in\DH$, such  that $L(S) = S\star T$ for all $S\in\Dip$. The eigenvalues are $m_\alpha=T_x\Big(\frac{\sigma(x)}{x^{\alpha+\be}}\Big)$.
\end{theorem}

\section{Extension of the operator $M$}

For $L\in\cM(\Om)$ we set $M=L^*\in L(\Do)$ and we study $L$ by means of properties of $M$. Due to reflexivity we have $L=M^*$. For $\vp\in\Do$ we set $\psi=M(\vp)$ and obtain the following characterization:
\begin{lemma}\label{l1} Let $L\in L(\Dop)$, then $L\in\cM(\Om)$ with eigenvalues $m_\alpha$ if and only if the following holds
\begin{equation}\label{eq1}
\int\xi^\alpha (m_\alpha \vp(\xi) - \psi(\xi)) d\xi =0
\end{equation}
for all $\alpha\in\N_0^d,\,\vp \in\Do$.
\end{lemma}

\Proof The condition is equivalent to
$$\int (m_\alpha\xi^\alpha) \vp(\xi) d\xi =\int \xi^\alpha (M\vp)(\xi) d\xi= \int L(\xi^\alpha) \vp(\xi) d\xi$$
for all $\vp\in\Do$. That is $m_\alpha \xi^\alpha = L(\xi^\alpha)$ for all $\alpha\in\N_0$. \qed

Let now the family $m_\alpha$, $\alpha\in \N_0^d$, of real numbers be given.

\begin{definition}\label{d5} We set
$$D(\tM)=\{\vp\in\Di\,:\,\ex \psi\in\Di\text{ such that } (\ref{eq1}) \text{ holds }\fa \alpha\in\N_0^d\}$$
and
$$\Gamma(\tM)=\{(\vp,\psi)\in\Di\times\Di\,:\,\vp,\,\psi\text{ fulfill }(\ref{eq1}) \fa \alpha\in\N_0^d\}.$$
\end{definition}

\begin{lemma}\label{l2} (1) $D(\tM)$ is a linear space, closed under dilations.\\
(2) $\Gamma(\tM)$ is the graph of a linear map $\tM:D(\tM)\to\Di$ which commutes with dilations.
\end{lemma}

\Proof We show only that $\Gamma(\tM)$ is a graph, the rest is obvious. Because of linearity of $\Gamma(\tM)$ it suffices to assume $\vp=0$. Then (\ref{eq1}) says that $\int\xi^\alpha \psi(\xi) d\xi=0$ for all $\alpha\in\N_0^d$, hence $\psi=0$. \qed

For open $\Om$ we set $\tOm=\bigcup_{a\in\NZ}a\Om$. By definition $\tOm$ is invariant under dilations. We obtain:

\begin{lemma}\label{l6} 1. If $\Do\subset D(\tM)$ then $\cD(\tOm)\subset D(\tM)$.\\
 2. If $\tM\cD(\Om)\subset \cD(\Om)$ then $\tM\cD(\tOm)\subset \cD(\tOm)$ and $\tM \in L(\cD\,'(\tOm))$.
\end{lemma}

\Proof By Lemma \ref{l2} we have $\cD(a\Om)\subset D(\tM)$ for all $a\in\NZ$. If $\vp\in\cD(\tOm)$ then there are finitely many $a_j\in\NZ$ such that $\vp\in\cD(\bigcup_j a_j\Om)$. Hence there are $\vp_j\in\cD(a_j\Om)$ such that $\vp=\sum_j\vp_j$. Since $D(\tM)$ is a linear space 1. is proved. 2. follows from the fact that $\tM (a\Om)\subset a\Om$ for all $a\in\NZ$ and the above decomposition, since $\tM$ is a linear map. The continuity of $\tM:\cD(\tOm)\to \cD(\tOm)$ follows from de Wilde's Theorem (or Grothendieck's closed graph theorem). \qed

\begin{lemma}\label{l4} $\NZ\subset \tOm$ for all non-empty open $\Om\subset\Rn$.
\end{lemma}

\Proof Let $x_0\in\Om\cap\NZ$. Then $\NZ=\{ax_0\,:\,a\in\NZ\}\subset\tOm$. \qed

We will later study the following two special cases:

\begin{remark}\label{r2} 1. If\, $0\in\Om$ then $\tOm=\Rn$.\\ 2. If $\Om\subset \NZ$ then $\tOm=\NZ$.
\end{remark}

\section{Representation}

In this whole section we assume that $L\in\cM(\Om)$, $M=L^*$ and $m_\alpha$, $\alpha\in\N_0^d$, the family of eigenvalues. In this case $\Do\subset D(\tM)$ and $M\vp=\tM\vp$ for $\vp\in\Do)$.

 We will use the following notations:

For $\eps>0$ we set $$W_\eps=\{x\in\Rn\,:\, \min_j |x_j|\ge\eps\}.$$

  For $r=(r_1,\dots,r_d)$, where all $r_j\ge 0$, we set $B_r=\{x\in\R^d:|x_j|\le r_j\text{ for all }j\}$ and for $s\ge 0$ we set $B_s:=B_{s\be}=\{x\in\Rn\,:\,|x|_\infty\le s\}$.\\
For $r=(r_1,\dots,r_d)\in\Rn$ and $s\in\R$ we set $r+s=r+s\be=(r_1+s,\dots,r_d+s)$.

 Since $\tM$ commutes with dilations we have
 $$\tM_\xi(\vp(\eta\xi))[x]=(\tM\vp)(\eta x)$$
 for all $\vp\in\cD(\tOm)$ and $\eta\in\NZ$.

 For $\vp\in\cD(\tOm)$ we define
 $$T\vp = (\tM\vp)(\be).$$
 Then $T\in\cD\,'(\tOm)$ and for all $\eta\in\NZ$ and $\vp\in\cD(\tOm)$ we have $\vp(\eta\,\cdot) \in\cD(\tOm)$ and
 \begin{equation}\label{eq2} (\tM\vp)(\eta)=T_\xi \vp(\eta\xi).
 \end{equation}

 We will study the properties of $T$. First we will investigate the consequences of the fact that $\eta\mapsto T_\xi \vp(\eta\xi)\in\Di$ for all $\vp\in\cD(\tOm)$. For that we modify the proof of Lemma 2.1 in \cite{Vdprime}.

 \begin{lemma}\label{l10} If $T\in\cD\,'(\tOm)$ and $\eta\mapsto T_\xi \vp(\eta\xi)\in\Di$ for all $\vp\in\cD(\tOm)$ then $\supp T\subset W_\eps$ for some $\eps>0$.
 \end{lemma}

 \Proof
 For $s\in\R_+^d$ we set $$K_s=\prod_j[\frac{s_j}{2},s_j]$$
 and we set $K_1=K_{(1,\dots,1)}$.

Since, by Lemma \ref{l6}, the map $\tM$ is continuous there is $r>0$ such that $\tM(\cD(K_1))\subset \cD(B_r)$ and this implies that $T_\xi\vp(\eta\xi)=0$ for any $\vp\in\cD(K_1), \eta\in\NZ$ and $|\eta|_\infty>r$.

We set $\eps=1/r$ and assume that $K_s\subset \{x\,:\, x_j<\eps\}$. Then $|1/s|_\infty>r$. For $\vp\in\cD(K_s)$ we set $\psi(\xi)=\vp(\xi s)$. Then $\psi\in\cD(K_1)$ and
$$T\vp =T_x(\psi(x/s))=0.$$
We have shown that $T|_{\{x\in\R_+^d\,:\,x_j<\eps\}}=0$. Repeating this in an analogous way for all `quadrants' and all relevant half-spaces, we obtain $T|_{\NZ\setminus W_\eps}=0$.

For the support of $T$ as a distribution in $\cD(\tOm)$ we obtain
$$\supp T \subset (W_\eps\cap\tOm)\cup (\Rn\setminus\NZ)\cap\tOm.$$
So $T$ decomposes into two distribution $T_1$ with support in $W_\eps$ and $T_2$ with support in $\tOm\setminus\NZ$. It is well known (see \cite{LSI}) that $T_2$ is locally a finite sum of distributions of the form
$$\vp\mapsto S_{\xi'}\frac{\partial^\nu}{\partial \xi_j^\nu}\vp(\xi_1,\dots,\xi_{j-1},0,\xi_{j+1},\dots,\xi_d)$$
where $\xi'=(\xi_1,\dots,\xi_{j-1},\xi_{j+1},\dots,\xi_d)$ and $S$ a distribution in these variables. This implies that  $S_x\vp(xy)$ does not depend on $y_j$. Therefore $T_2$ must be $0$. \qed

We have shown:
\begin{proposition}\label{p2} There is a distribution $T\in\cD\,'(\tOm)$ with support in $W_\eps$ for some $\eps>0$ such that $(\tM\vp)(y)=T_x\vp(xy)$ for all $\vp\in\cD(\tOm)$ and $y\in\NZ$. This holds, in particular, for $\vp\in\cD(\NZ)$ and for $\vp\in\cD(\Om)$. In the latter case this means $(M\vp)(y)=T_x\vp(xy)$.
\end{proposition}

Next we want to get information on the support of $T$.  We need a preparatory lemma.

We will use the following notation: For $M,N\subset\Rn$ we set $M^c:=\Rn\setminus M$ and $V_*(M,N)=\{\eta\in\NZ\,:\,\eta M\subset N\}$, $V_*(M)=V_*(M,M)$.

\begin{lemma}\label{r3} Let $M,N\subset\Rn$ then $V_*(M^c,N^c)=\{1/y\,:\,y\in V_*(N,M)\}$.
\end{lemma}

\Proof The statement is symmetric, hence we need only to show one implication. Let $y\in V_*(N,M)$ then $M=y N\,\dot{\cup}\,y N_0$ where $N_0\subset N^c$, hence $\frac{1}{y}\,M=N\,\dot{\cup}\,N_0$, and this implies $\frac{1}{y}\,M^c = (N\,\dot{\cup}\,N_0)^c\subset N^c$. \qed

The following example shows that for the sets $V(M,N)=\{\eta\in\Rn\,:\,\eta M\subset N\}$ there is no such simple relation. For our theory we are only interested in $V_*(M,N)$.

\begin{example}\label{ex1} Let $\Om=\Rn\setminus \{0\}$ then $V(\Om)=\NZ$ while $V(\Om^c)=V(\{0\})=\Rn$.
\end{example}

For $M\subset\NZ$ we set $1/M=\{1/y\,:\, y\in M\}$.

\begin{proposition}\label{l7} If $T\in\cD\,'(\tOm)$ such that $T_x\vp(x \cdot)\in\Di$ for all $\vp\in \cD(\tOm)$, then the following are equivalent:
\begin{enumerate}
\item $T_x\vp(x \cdot)\in\Do$ for all $\vp\in \Do$.
\item For every $\om\ssubset\Om$ there is $L\ssubset\Om$ such that $\supp T\subset V_*(L^c,\om^c)$.
\item $(1/\supp T) K\ssubset\Om$ for every compact $K\subset \Om$.
\end{enumerate}
\end{proposition}

\Proof We first remark that, by Lemma \ref{l10} and the assumption, $\supp T\subset\NZ$.

1.$\Rightarrow$ 2.\; Let  $\om$ be an open set such that $\om\ssubset\Om$. Then there is a compact set $L\subset\Om$ such that $T_x\vp(xy)=0$ for all $\vp\in\cD(\om)$ and $y\not\in L$. If $y\in\NZ$ then
$$\{\vp(\cdot\,y)\,:\,\vp\in\cD(\om)\}=\cD\Big(\frac{1}{y}\om\Big).$$
Therefore
$$\supp T\cap\frac{1}{y}\om=\emptyset \text{ for all }y\in\NZ, y\not\in L.$$
 We have shown
$$\supp T \subset \bigcap_{y\in \NZ\cap L^c}\frac{1}{y}\, \om^c= V(L^c,\om^c).$$
The last equality comes from \cite[Lemma 1.4]{Vconv}. Since $\supp T\subset\NZ$ we can replace $V(L^c,\om^c)$ with $V_*(L^c,\om^c)$.

2.$\Rightarrow$ 1.\; Let $\om\subset\Om$ be compact and $\vp\in\cD(\om)$. We choose $L$ according to the assumption and assume that $\supp T\subset V_*(L^c,\om^c)$. If $x\in V_*(L^c,\om^c)$ and $y\in L^c$ then $xy\in \om^c$. Since $x\mapsto xy$ is continuous and $\om^c$ open we have $xy\in\om^c$, hence $\vp(xy)=0$, for $x$ in a neighborhood of $\supp T$. This shows that $T_x\vp(xy)=0$ for all $y\in L^c$.

2.$\Leftrightarrow$ 3.\, In 2. we might assume $\om$ and $L$ to be compact. So by Lemma \ref{r3} it is equivalent to:
For every compact $K\subset\Om$ there is a compact $L\subset\Om$ such that $1/\sup M\subset V_*(K,L)$. The equivalence to 3. is now obvious. \qed


From Proposition \ref{l7}, (3.) we get the following:
\begin{corollary}\label{c1} $1/\supp T\subset V_*(\Om)$.
\end{corollary}

In Proposition \ref{l7}, we cannot replace 3. with $1/\supp T\subset V_*(\Om)$ as the following example shows. See, however, Proposition \ref{l9}.

\begin{example}\label{ex4} Let $\Om=]0,1[$ and $T\vp = \int_1^\infty \vp(x) e^{-x} dx$. Then $\vp\mapsto T_x\vp(x\,\cdot)$ maps $\Di$ into $\Di$ (cf. \cite[Proposition 2.5]{Vdprime}). $1/\supp T =]0,1]=V_*(\Om)$. For $\vp\in\Do$, $\vp\ge 0$, $\vp\neq 0$  we obtain
$$T_x\vp(xy) = \frac{1}{y} \int_y \vp(x) e^{-x/y} dx$$
which is $>0$ near $0$, hence $T_x\vp(x\,\cdot)\not\in\Do$.
\end{example}

An important special case is the following.

\begin{theorem}\label{t9} If $\Om\ssubset Q_+$ then all Hadamard operators on $\Dop$ are Euler operators.
\end{theorem}

\Proof Let $L\in\cM(\Om)$, then by Proposition \ref{p2} and Corollary \ref{c1} there is $T\in\Dip$ with $1/\supp T\subset V_*(\Om)$ such that $(M\vp)(y)=T_x\vp(xy)$.  Since $V(\Om)=\{\be\}$ (see e.g. \cite[Example 5.3]{Vhad}) we have $\supp T =\{\be\}$ or $T=0$.
Therefore $T=\sum_\alpha b_\alpha \delta_\be^{(\alpha)}$ hence $(M\vp)(y)=\sum_\alpha b_\alpha (-1)^{|\alpha|} y^\alpha \vp^{(\alpha)}(y)$ and therefore $(LS)(\vp)=S(M\vp)=\sum_\alpha b_\alpha (-1)^{|\alpha|} S_y(y^\alpha \vp^{(\alpha)}(y))$ for all $\vp\in\Do$. This means that $L(S)=(\sum_\alpha c_\alpha \theta^\alpha) S$ with suitable constants $c_\alpha$. \qed

\begin{definition}\label{d6} For open $\Om\subset\Rn$ we define: $T\in\DHO$ if $T\in\cD\,'(\tOm)$ and
\begin{enumerate}
\item $T_x\vp(x,\,\cdot)\in \cD(\tOm)$ for all $\vp\in\cD(\tOm)$,
\item $(1/\supp T) K\ssubset\Om$ for every compact $K\subset \Om$.
\end{enumerate}
\end{definition}

We have shown one of our main results:

\begin{theorem}\label{t4} If $L\in\cM(\Om)$ then there is $T\in\DHO$ such that
$L(S)=S\star T$ for all $S\in\Dop$.
\end{theorem}

It remains to give a closer description of the distributions $T$ appearing in Definition \ref{d6}. This can be done it in the two important cases mentioned in  Remark \ref{r2}. In the case of $0\in\Om$ it follows from \cite{Vdprime}.

\section{Case of $0\in\Om$}

Let $\Om\subset\Rn$ be open with $0\in\Om$, then $\tOm=\Rn$ (see Remark \ref{r2}). Let $L\in\cM(\Om)$.  By Theorem \ref{t4} there is $T\in\DHO$ such that
$L(S)=S\star T$ for all $S\in\Dop$. Condition 1. in Definition \ref{d6} then means that $T\in\DH$, that is, $T\in\OH$ and $\supp T\subset W_\eps$ for some $\eps>0$. We have proved one implication of the following theorem.

\begin{theorem}\label{t5} If $0\in\Om$, then
$L\in\cM(\Om)$ if and only if
there is a distribution $T\in\DH$ with $(1/\supp T) K\ssubset\Om$ for every compact $K\subset \Om$, such that $L(S)=S\star T$ for all $S\in\Dop$.
\end{theorem}

\Proof It remains to show that under the given conditions $M_T:\vp\mapsto T_x\vp(x\,\cdot)$ is the transpose of an operator in $\cM(\Om)$. Clearly $M_T\in L(\Do)$. We have to verify the condition in Lemma \ref{l1} and this follows from the proof of \cite[Theorem 4.2]{Vdprime}. \qed

For the open unit ball and similar sets the we obtain (cf. \cite[Theorem 4.4]{Vhad}).

\begin{example}\label{ex2} Let $\Om\subset\Rn$ be a bounded open set  with the following properties:\\
1. If $x\in\Om$ and $|y_j|\le|x_j|$ for all $j$ then also $y\in\Om$.\\
2. $\Om$ is invariant under permutations of the variables.\\
3. $t\,\overline{\Om}\subset \Om$ for all $0<t<1$.\\
Then $\DHO=\{T\in\OH\,:\,\supp T\subset W_1\}$.
\end{example}

\Proof We may choose $\om=\eps\Om$, $L=\delta\Om$. Then
$$V(L^c,\om^c)=V(\delta \Om^c,\eps \Om^c)=\frac{\eps}{\delta} V(\Om^c).$$
By \cite[Theorem 4.4]{Vhad} we have $V(\Om)=Q:=\{x\,:\,|x|_\infty\le1\}$ hence $V_*(\Om^c)=1/V_*(\Om)= W_1$. This implies
$$V_*(L^c,\om^c) =\frac{\eps}{\delta} V_*(\Om^c)= \frac{\eps}{\delta} W_1.$$
So we obtain $T\in\DHO$ if and only if for every $0<\eps<1$ there is $0<\delta<1$ such that $\supp T\subset (\eps/\delta) W_1$ and this is the case if and only if $\supp T\subset W_1$. \qed

\begin{theorem}\label{t6} If $\Om$ is the open unit ball for $\ell_p$, $0<p\le+\infty$, then
$L\in\cM(\Om)$ if and only if
there is a distribution $T\in\OH$, $\supp T\subset W_1$, such that $L(S)=S\star T$ for all $S\in\Dop$.
\end{theorem}

\section{Case of $\Om\subset\NZ$}

Let $\Om\subset\NZ$ be open, then $\tOm=\NZ$ (see Remark \ref{r2}). Let $L\in\cM(\Om)$.  By Theorem \ref{t4} there is $T\in\DHO$ such that
$L(S)=S\star T$ for all $S\in\Dop$. In particular, $T\in\cD\,'_H(\NZ)$.

\begin{lemma}\label{l8} $\cD\,'_H(\NZ) = \cE'(\NZ)$
\end{lemma}

\Proof Assume that $T\in\cD\,'_H(\NZ)$ that is $\supp T\subset W_\eps$ for some $\eps>0$ and $(\tM\vp)(y)=T_x\vp(xy)$ for all $\vp\in\cD(\NZ)$ and $y\in\NZ$.

We set $\check{\vp}(x)=\vp(1/x)$ and $\check{T}\vp = T\check{\vp}$ for $\vp\in\cD(\NZ)$. Then $\check{T}\in \cD\,'(\NZ)$. We obtain for $\vp\in\cD(\NZ)$ and $y\in\NZ$
$$\check{T}_x\vp(xy)=T_x\vp\Big(\frac{1}{x}\,y\Big)=T_x\check{\vp}\Big(x\,\frac{1}{y}\Big).$$
$T_x\check{\vp}(x\,\frac{1}{y})\in\cD(\NZ)$ as a function of $y$. By Proposition \ref{p2} there is $\delta>0$ such that $\supp \check{T}\subset W_\delta$.

We set $r=1/\delta$ and $B_r=\{x\,:\,|x|_\infty \le r\}$. Assume that $\supp \vp\in\NZ\setminus B_r$, then $\supp\check{\vp}\in\NZ\setminus W_\delta$. Therefore $T\vp=\check{T}\check{\vp}=0$.

We have shown that $\supp T \subset W_\eps\cap B_r$, hence compact. The other implication is obvious.   \qed

We refer to Corollary \ref{c1} and show for $\Om\subset\NZ$:

\begin{proposition}\label{l9} $T\in\DHO$ if and only if $T\in\cE'(\NZ)$ and $1/\supp T\subset V_*(\Om)$.
\end{proposition}

\Proof The `only if' part follows from the Definition \ref{d6} and Corollary \ref{c1}. For the reverse direction we note that, by Lemma \ref{l8}, $T\in\cD\,'_H(\NZ)$.
Since $\supp T\subset \NZ$ is compact, also $1/\supp T$ is compact. Therefore $(1/\supp T) K$ is compact for any compact $K\subset\Om$ and, by assumption, it is contained in $\Om$. Therefore $T\in\DHO$. \qed

So we have shown one implication of the following theorem.

\begin{theorem}\label{t7} If $\Om\subset\NZ$, then
$L\in\cM(\Om)$ if and only if
there is a distribution $T\in\cE'(\NZ)$ with  $1/\supp T\subset V_*(\Om)$, such that $L(S)=S\star T$ for all $S\in\Dop$.
\end{theorem}

\Proof It remains to show that under the given conditions $M_T:\vp\mapsto T_x\vp(x\,\cdot)$ is the transpose of an operator in $\cM(\Om)$. Clearly $M_T\in L(\Do)$. We have to verify the condition in Lemma \ref{l1} and this follows from the proof of \cite[Theorem 4.1]{Vdprime}. \qed

\section{Final remarks}

If $d=1$ the the cases in Theorems \ref{t5} and \ref{t7} cover all possibilities, hence we have a complete characterization.

\begin{theorem}\label{t8} If $\Om\subset \R$ is open then
$L\in\cM(\Om)$ if and only if
there is a distribution\\
$T\in\Dh$ with $(1/\supp T) K\ssubset\Om$ for every compact $K\subset \Om$\\ if $0\in\Om$
or\\
 $T\in\cE'(\NZ)$ with $\supp T\subset \NZ$ and $1/\supp T\subset V_*(\Om)$\\ if $\Om\subset \R_*$\\
 such that
$L(S)=S\star T$ for all $S\in\Dop$.
\end{theorem}

An interesting example in higher dimensions is $\Om=\Rn\setminus\{0\}$:

\begin{example}\label{ex3}  In this case we have also
$\cD\,'_H(\Rn\setminus\{0\}) = \cE'(\NZ)$ and the proof is like in the proof of Lemma \ref{l8} except we habe to replace the pointwise reciprocal by the reflection at the euclidian unit sphere $x\mapsto x/|x|^2$.
\end{example}

If, for instance, $L\in\cM(\Om)$, $\Om\subset\R^2$ and $0\not\in\Om$, $e_1,e_2\in\Om$ where $e_j$ are the unit vectors, then $\tOm=\R^2\setminus \{0\}$ and we have for the representing distribution $T$ the condition $T\in\cE'(\NZ)$ plus the support condition.

\vspace{.5cm}

\noindent Bergische Universit\"{a}t Wuppertal,
\newline Fak. Math.-Nat., Gau\ss -Str. 20,
\newline D-42119 Wuppertal, Germany
\newline e-mail: dvogt@math.uni-wuppertal.de


\begin{thebibliography}{99}

 \bibitem{DLre} P. Doma{\'n}ski, M. Langenbruch, Representation of multipliers on spaces of real
analytic functions, {\em Analysis} {\bf 32} (2012), 137--162.

\bibitem{DLalg} P. Doma{\'n}ski, M. Langenbruch, Algebra of multipliers on the space of real
analytic functions of one variable, {\em Studia Math.} {\bf 212} (2012), 155--171.

\bibitem{DLhad} P. Doma{\'n}ski, M. Langenbruch, Hadamard multipliers
on spaces of real analytic functions, {\em Adv. Math.} {\bf 240}
(2013), 575--612.

\bibitem{DLV} P. Doma{\'n}ski, M. Langenbruch, D. Vogt,
 Hadamard type operators  on spaces of real analytic functions in several variables, {\em J. Funct. Anal.} {\bf 269} (2015),  3868-–3913.

 \bibitem{DK} J. J. Duistermaat, J. A. C. Kolk, {\em Distributions. Theory and Applications}, Birkh\"auser, Boston 2010

 \bibitem{MV} R. Meise, D. Vogt: \it Introduction to functional analysis, \rm Clarendon Press, Oxford, (1997).

  \bibitem{LSI} L. Schwartz, \em Th\'eorie des distributions I, \rm Hermann, Paris 1957.

\bibitem{LS} L. Schwartz, \em Th\'eorie des distributions II, \rm Hermann, Paris 1959.

 \bibitem{Vhad} D. Vogt, Operators of Hadamard type on spaces of smooth functions, {\em Math. Nachr.}, {\bf 288} (2015), 353--361.

\bibitem{Vconv} D. Vogt, $\cE'$ as an algebra by multiplicative convolution, {\em arXiv: 1509.05759} (2015).

\bibitem{Vdprime} D. Vogt, Hadamard operators on $\cD\,'(\Rn)$, {\em arXiv: 1511.08593} (2015).

\end{thebibliography}
\end{document}